\documentclass[12pt,a4paper]{amsart}
\usepackage[latin1]{inputenc}
\usepackage{mathptmx}
\usepackage{amssymb}
\usepackage{hyperref}
\usepackage{breakurl}
\usepackage{mathrsfs}
\usepackage{graphicx}
\usepackage[mathscr]{euscript}
\usepackage{booktabs}

\usepackage{pgf, tikz}

\usepackage{fancyvrb} \RecustomVerbatimEnvironment{Verbatim}{Verbatim}
{xleftmargin=15pt, frame=single, fontsize=\small}

\usepackage{color}

\newtheorem{theorem}{Theorem}

\newtheorem{proposition}[theorem]{Proposition}
\newtheorem{conjecture}[theorem]{Conjecture}

\theoremstyle{definition}
\newtheorem{remark}[theorem]{Remark}

\numberwithin{equation}{section}

\def\RR{{\mathbb R}}

\def\11{{\mathbb 1}}

\def\cE{{\mathscr E}}

\def\Sg{{\textup {Sg}}}
\def\SgRev{{\textup {SgRev}}}

\def\CW{\textup{CW}}

\def\vol{\operatorname{vol}}

\def\prob{\operatorname{prob}}

\def\ttt#1{\texttt{#1}}

\let\epsilon=\varepsilon

\allowdisplaybreaks

\begin{document}





\begin{center}
\Large\textbf{Computations of volumes in five candidates elections$^{\dagger}$}\\
\bigskip
\large{Winfried Bruns$^{1,*}$ and Bogdan Ichim$^{2,3}$}\\
\bigskip
\end{center}
$^{1}$Universit\"at Osnabr\"uck, Institut f\"ur Mathematik, 49069 Osnabr\"uck, Germany, email: wbruns@uos.de.\\
$^{2}$University of Bucharest, Faculty of Mathematics and Computer Science, Str. Academiei 14, 010014 Bucharest, Romania, email: bogdan.ichim@fmi.unibuc.ro.\\
$^{3}$ Simion Stoilow Institute of Mathematics of the Romanian Academy, Research Unit 2, C.P. 1-764, 010702 Bucharest, Romania, email: bogdan.ichim@imar.ro.\\
$^{*}$Corresponding Author.\\
$^{\dagger}$The article has appeared in \emph{Scientific reports} 
\url{https://www.nature.com/articles/s41598-023-39656-8}. There are typographical differences between the two versions.
\bigskip

\begin{center}
	\emph{To the memory of Udo Vetter, our teacher, colleague and friend}
\end{center}
\bigskip

\textbf{Abstract.} We describe several analytical (i.e., precise) results obtained in five candidates
social choice elections under the assumption of the Impartial Anonymous Culture. These include the Condorcet and Borda paradoxes, as well as the Condorcet efficiency of plurality, negative plurality and Borda voting, including their runoff versions. The computations are done by Normaliz. It finds precise probabilities as volumes of polytopes in dimension~$119$, using its recent implementation of the Lawrence algorithm.
\bigskip

\section*{Introduction}

In \cite[p.\ 382]{LLS} Lepelley, Louichi and Smaoui state:
\begin{quotation}
	``Consequently, it is not possible to analyze
	four candidate elections, where the total number of variables (possible
	preference rankings) is $24$. We hope that further developments of these algorithms
	will enable the overcoming of this difficulty.''
\end{quotation}
This hope has been fulfilled by previous versions of Normaliz \cite{Nmz}. In connection with the symmetrization suggested by Sch\"urmann \cite{Sch}, it was possible to compute volumes and Ehrhart series for many voting events in four candidates elections; see \cite{BIS4cand}. As far as Ehrhart series are concerned, we cannot yet offer progress. But the volume computation was already substantially improved by the descent algorithm described in \cite{Desc}.
Examples of Normaliz being used for voting theory computations by independent authors can be found in \cite{BGS}, \cite{BHS} and \cite{Diss}.
The purpose of this paper is to present \emph{precise} probability computations in five candidates elections under the assumption of the Impartial Anonymous Culture (IAC). They are made possible by Normaliz' implementation  of the Lawrence algorithm \cite{Law}.

The connection between rational polytopes and
social choice was established independently in \cite{LLS} and  \cite{WP}.
Solutions for the four candidates quest were proposed for example in  \cite{Sch}, \cite{BIS4cand} and \cite{Desc}.
The similar, but much more challenging computational problem of performing precise computations in five candidates elections
is wide open. Various authors have used the well known Monte Carlo methods in order to perform computations with five or more candidates,
but fundamentally these methods can only deliver \emph{approximative} results, without even clear bounds for errors. We note that methods that were successful in obtaining precise results
in the four candidates case are ineffective in the five candidates case due to the huge leap in computational complexity implied by the increase in the dimension of the associated polytopes (from $23$ to $119$).
Therefore a different algorithmic approach is needed in order to obtain the desired precise results.

To the best of our knowledge, we present here the first precise results obtained for computations with five candidates.
By precise we mean either absolutely precise rational numbers, or results obtained using the fixed precision mode of Normaliz where the desired precision is set and fully controlled by the user.

The polytopes in five candidates elections have dimension $119$, and are defined as subpolytopes of the simplex spanned by the unit vectors of $\RR^{120}$. The number of the inequalities cutting out the subpolytope is the critical size parameter, but fortunately we could manage computations with $\le 8$ inequalities (in addition to the $120$ sign inequalities) on the hardware at our
disposal, although the algorithm allows an arbitrary number of inequalities. This covers the Condorcet paradox \cite{C} (computable on a laptop in a few minutes), the Borda winner and loser paradoxes \cite{B}, and the Condorcet efficiency of plurality, negative plurality and Borda voting, including their runoff extensions. We also compute the probabilities of all $12$ configurations of the five candidates that are defined by the Condorcet majority relation.

As Table \ref{calculation_times} shows, the computations for $5$ candidates are very demanding on the hardware in memory and computation time. Therefore we consider it a major value of the new algorithm that it improves the situation in four candidates elections considerably, where it is now possible to allow preference rankings with all types of partial indifference. Moreover one can run series of parameterized computations for four candidates like those that one finds in \cite{Indiff} for three candidates. In order to illustrate this possibility we compute the probability of the Condorcet paradox in the presence of voters with indifference and the Condorcet efficiency of approval voting (see Subsection "Indifference").
Note that potential applications are not only limited to voting theory, as can be seen in \cite[Table 3]{IM}. There the new algorithm is performing better (as the dimension grows) for the first family of examples.

Normaliz computes lattice normalized volume and uses only rational arithmetic without rounding errors or numerical instability. But there is a slight restriction: while it is always theoretically possible to compute the probabilities as absolutely precise rational numbers, the fractions involved can reach sizes which are unmanageable on the available hardware. For these cases Normaliz offers a fixed precision mode whose results are precise up to an error with a controlled bound that can be set by the user.

In contrast to algorithms that are based on explicit or implicit triangulations of the polytope $P$ (or the cone $C(P)$ defined by $P$) under consideration, the Lawrence algorithm uses a ``generic triangulation'' of the dual cone $C(P)^*$. We make a brief discussion of the available Lawrence algorithm implementations and their limitations in Section "Implementations of the Lawrence algorithm and their limitations". In order to reach the order of magnitude that is necessary for five candidates elections, one needs a fine tuned implementation. It is outlined in \cite{BVol}. Moreover, the largest of our computations need a high performance cluster to finish in acceptable time. Section "Computational report" gives an impression on the computation times and memory requirements by listing them for selected examples.

The computations that we report in this note were done by version 3.9.0 of Normaliz. Meanwhile it has been succeeded by version 3.10.1 without changes in the Lawrence algorithm. Both versions are available at
\begin{center}
\url{https://www.normaliz.uni-osnabrueck.de/}
\end{center}
For details on the implementation and the performance of the previous versions of Normaliz we point the reader to \cite{BK}, \cite{BI}, \cite{BIS}, \cite{BS}.

\section*{A Challenging Computational Problem Arising from Social Choice}

\subsection*{Voting schemes and rational polytopes}\label{scheme_vol}

The connection between voting schemes and rational polytopes is based on counting integral points in the latter. In this subsection we sketch the connection. As a general reference for discrete convex geometry we recommend \cite{BG}. The interested reader may also consult \cite{GL} and \cite{GLnew}.

The basic assumption in the mathematics of social choice is the existence of \emph{individual preference rankings} $\succ$: every voter ranks the candidates in linear order. Examples for three candidates named by capital letters:
$$
A \succ B \succ C,\qquad C\succ A\succ B.
$$
For $n$ candidates there exist $N=n!$ preference rankings, usually numbered in lexicographic order.  (By an extension it is possible to allow indifferences; for example see \cite{Indiff}.)

The \emph{result} or \emph{profile} of the election is the $N$-tuple
$$
(x_1,\dots,x_N),\qquad x_i=\#\{  \textup{voters of preference ranking } i\}.
$$
Thus an election result for three candidates may be written in the following tabular form:
\begin{center}
	\begin{tabular}{c||c|c|c|c|c|c}
		number of voters&$x_{1}$&$x_{2}$&$x_{3}$&$x_{4}$&$x_{5}$&$x_{6}$\\
		\hline
		&$A$&$A$&$B$&$B$&$C$&$C$\\
		ranking&$B$&$C$&$A$&$C$&$A$&$B$\\
		&$C$&$B$&$C$&$A$&$B$&$A$
	\end{tabular}
\end{center}

In the following we want to compute probabilities of certain events related to election schemes. This requires a probability distribution on the set of election results. The \emph{Impartial Anonymous Culture (IAC)} assumes that all election results  for a fixed number of voters, in the following denoted by $k$, have equal probability. In other words, it is the equidistribution on the set of voting profiles for a fixed number of $k$ voters.

The \emph{Marquis de Condorcet} (1743--1794) was a leading intellectual in France before and during the revolution. He already observed that there is no ideal election scheme, a fact now most distinctly manifested by Arrow's impossibility theorem.
We say that candidate $A$ \emph{beats} candidate $B$ in \emph{majority}, $A>_M B$,  if
$$
\#\{\textup{voters with } A\succ B \} > \#\{\textup{voters with } B\succ A \}.
$$
A (necessarily unique) \emph{Condorcet winner} (CW) beats all other candidates in majority. There is general agreement that the CW is the person with the largest common approval.
However, Condorcet realized that a CW need not exist: the relation $>_M$ is not transitive: a minimal example is the profile $(1,0,0,1,1,0)$. This phenomenon is called the \emph{Condorcet paradox}. From a quantitative viewpoint, the most ambitious goal is to find the exact number of election profiles exhibiting the Condorcet paradox (or the opposite), given the number of voters $k$.
For  large $k$, this number is gigantic. It is much more informative to understand the behavior for $k\to\infty$:
what is the probability that an election result exhibits the Condorcet paradox?
Since we assume the IAC, this probability is
$$
\lim_{k\to\infty} \frac{\#\{\textup{electionresults without CW for $k$ voters}  \}}{\#\{\textup{all election results for $k$ voters}  \}}.
$$
It is a crucial consequence of (IAC) that the event ``$A$ is the CW'' can be characterized by a system of homogeneous linear inequalities. For three candidates they are
\begin{align*}
	A>_M B:\quad x_{1}+x_{2}+x_{5} &> x_{3}+x_{4}+x_{6},\\
	A>_M C:\quad x_{1}+x_{2}+x_{3} &>x_{4}+x_{5}+x_{6}.
\end{align*}
If we are only interested in probabilities for $k\to \infty$, standard arguments of measure theory  allow ties and replacement of $>$ by~$\ge$.

\begin{figure}[hbt]
\begin{center}
	\tikzset{facet style/.style={opacity=1.0,very thick,line,join=round}}
	\begin{tikzpicture}[x  = {(-0.5cm,-0.5cm)},
		y  = {(0.9659cm,-0.25882cm)},
		z  = {(0cm,1cm)},
		scale = 2]
		\draw [->] (-0.5, 0, 0) -- (1.7,0,0) node at (1.9,0,0) {$x_1$} node at (1, -0.1,0){$k$};
		\draw [->] (0, -0.5, 0) -- (0,1.5,0) node at (0,1.7,0) {$x_2$} node at (-0.1, 1,0.1){$k$};
		\draw [->] (0, 0, -0.5) -- (0,0,1.5) node at (0,0,1.6){$x_3$}node at (-0.1, -0.15,1){$k$};
		\draw[thick] (1,0,0) -- (0,1,0) -- (0,0,1) -- cycle;
		\filldraw[color=gray] (0.5,0.5,0) -- (0.266, 0.4, 0.333) -- (0,0.5,0.5)  -- (0,1,0) -- cycle;
		\draw node at (1,0,1){$\Delta_k$};
		\draw node at (0,0.8,0.4){$P_k$};
	\end{tikzpicture}
\end{center}
\caption{Subpolytope defined by linear inequalities}\label{subpolytope}
\end{figure}
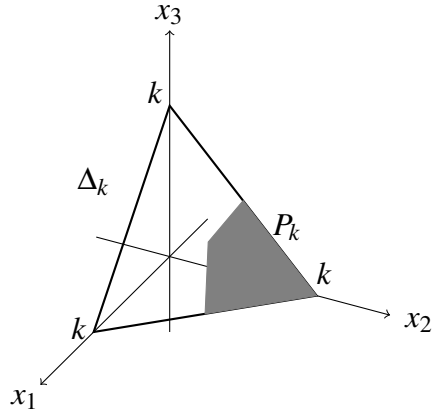

We now consider an event $E$ defined for an $n$ candidates election by a system of homogeneous linear inequalities on the set of election profiles. As above, set $N=n!$. The election profiles $(x_1,\dots,x_N)$ are the lattice points (points with integral coordinates) in the positive orthant $\RR_+^N$ satisfying the equation $x_1+\dots+x_N=k$. The real points in the positive orthant satisfying this equation form a polytope $\Delta_k$, and the linear inequalities whose validity defines $E$ cut out a subpolytope $P_k$. We illustrate this assertion by the (necessarily unrealistic) Figure \ref{subpolytope}.

For large numbers of voters we want to find the probability $\prob(E)$ of the event $E$. Under (IAC) it is given by
$$
\prob(E)=\lim_{k\to\infty}\frac{\#\{\textup{lattice points in } P_k \} }  {\#\{ \textup{lattice points in } \Delta_k \} }.
$$
We project $\Delta_k$ orthogonally onto $\Delta_1$, and thus $P_k$ onto $P_1$. The density, roughly speaking, of the projections of the lattice points converges to $1$, and therefore
$$
\prob(E) =\frac{\vol(P_1)}{\vol(\Delta_1)}.
$$
For volume computations in connection with the counting of lattice points one uses the \emph{lattice normalized volume} $\vol$, giving volume $1$ to $\Delta_1$. With this choice $\prob(E) = \vol(P_1)$.

It is not difficult, but would take many pages, to write down the linear inequalities for the voting schemes and events discussed in the following. For four candidates the complete systems are contained in \cite{BIS4cand}. For the inequalities one must often fix the roles that certain candidates play, like the Condorcet winner $A$ above. Then probabilities must be computed carefully, and this may require the inclusion-exclusion principle.

Both from the theoretical as well as from the computational viewpoint it is better to consider the cone $C$ defined by the homogeneous linear inequalities as the prime object, and the polytopes as intersections of $C$ with the hyperplane defined by the equation $x_1+\dots + x_N= k$.

It is not difficult to see that a voting event that can be realized by a voting profile has positive probability:

\begin{proposition}\label{pos_prob}
Let $\cE$ be a subset of all voting profiles defined by strict homogeneous rational inequalities. If $\cE$ is nonempty, then it has probability $> 0$ under (IAC).
\end{proposition}

\begin{proof}
Clearing denominators, one can assume that the coefficients of the inequalities are integers. Let $m$ be the maximum of all their absolute values and $x\in\cE$ be a voting profile. Then $x'=(m+1)x \in \cE$ as well by homogeneity. It is easily checked that also $x'+e_i\in\cE$ where $e_i$, $i=1,\dots,N$ is the $i$-th unit vector. The parallel translation by $-x'$ maps the the polytope $P$ spanned by the $x'+e_i$ bijectively onto $\Delta_1$. Thus $P$ has lattice normalized volume $1$, and therefore its orthogonal projection to $\Delta_1$ has positive volume.
\end{proof}

\subsection*{The Condorcet paradox in five candidates elections}  The Condorcet paradox, introduced in Subsection "Voting schemes and rational polytopes", does not occur in the case of two candidates
(if draws are excluded). For three candidates the exact probability of an outcome with a Condorcet winner (under IAC) was first computed by Gehrlein and Fishburn \cite{GF}
while for four candidates it was first determined by Gehrlein in \cite{G}.

For five candidates, we have computed in the full precision mode of Normaliz (and the method presented in Section "Implementations of the Lawrence algorithm and their limitations") that
$$
p_\CW=\frac{a}{b},
$$
where
$$
	a= \scalebox{0.9}{760794547958864241496408591531018198021484884229346111658236615929935}
$$
and
$$
	b= \scalebox{0.9}{1010827262551214358630401511004028249102084136257935356483840264634368.}
$$
In decimal notation with $100$ decimals, we obtain
\begin{align*}
	p_\CW \approx\, & \scalebox{0.9}{0.75264545797736434427639219331756247271265813365410} \\
	& \scalebox{0.9}{18228684464583400970327543361542592465345709195008.}
\end{align*}

In order to illustrate the fixed precision mode of Normaliz, we compare the above exact result with the result obtained for fixed precision of $100$ decimal digits, namely
$$
p_\CW\approx\frac{a'}{b'},
$$
where
\begin{align*}
	a'=\, & \scalebox{0.9}{940806822471705430345490241646953090890822667067627278585558072925121}\\
    & \scalebox{0.9}{29094292019282405816821364997}
\end{align*}
and
\begin{align*}
	b'=\, & \scalebox{0.9}{125000000000000000000000000000000000000000000000000000000000000000000}\\
    & \scalebox{0.9}{000000000000000000000000000000.}
\end{align*}
In decimal notation with $100$ decimals, we obtain
\begin{align*}
	p_\CW \approx & \scalebox{0.9}{0.75264545797736434427639219331756247271265813365410} \\
	& \scalebox{0.9}{18228684464583400970327543361542592465345709199760.}
\end{align*}
The reader should observe that in the decimal notation only the last $4$ digits are different. The error bound is
$$
6572904\cdot 10^{-100} < 10^{-93},
$$
where $6,572,904$ is the size of the ``generic triangulation'' (see Section "Implementations of the Lawrence algorithm and their limitations" and Table  \ref{numerical_data}).

This means that using the fixed precision mode of Normaliz is sufficient for many applications,
while it saves computation time and is significantly less demanding on the hardware.
\bigskip

For practical reasons, in the following we use shorter decimal representations of the rational numbers. (The full rational representations of these numbers are available on demand from the authors.) A decimal representation is called \emph{rounded to $n$ decimals} when the first $n-1$ printed decimals are exact and only the last decimal may be rounded up.

\subsection*{Rule versus rule runoff, Condorcet efficiencies} The most common voting scheme in elections is the \emph{plurality rule} $PR$: for each candidate $X$ one counts the voters that have $X$ on first place in their preference ranking, and the winner is the candidate with most first places. However, in many elections one uses a second ballot, called \emph{runoff}, if the winner has not got the votes of more than half of the voters. In the runoff only two candidates are left, namely the two top candidates of the first round. A typical example is the French presidential election.

If the ideal winner of an election is the Condorcet winner CW,  then one must ask for the probability that the plurality winner is the CW under the condition that a CW exists. This conditional probability is called the \emph{Condorcet efficiency}, studied intensively by Gehrlein and Lepelley \cite{GL} as a quality measure for voting schemes.

Another important question is whether the runoff is a real improvement: (i) what is the probability that the winner of the first ballot also wins the second, and (ii) by how much does the Condorcet efficiency increase by the runoff.

An often discussed variant of plurality is \emph{negative plurality} $NPR$: the winner is the least disliked candidate $X$, defined by the least number of voters who have placed $X$ on the last place in their preference ranking. As for plurality one can have a runoff, and again it makes sense to compute the Conndorcet efficiencies and the probability that the first round winner also wins the runoff.

Both plurality and negative plurality are special cases of weighted voting schemes in which the places in the preference ranking have a fixed weight, and every candidate is counted with the sum of the weights in the preference ranking of the voters. In plurality the first place has weight $1$ and the other places have weight $0$, wheres negative plurality gives weight $-1$ to the last place. In addition to these two rules we discuss the \emph{Borda rule} $BR$ that for $n$ candidates gives weight $n-p$ to place $p$.

In the case of four candidates the plurality voting versus plurality runoff problem was first computed by De Loera, Dutra, K\"oppe, Moreinis, Pinto and Wu in \cite{Latte} using LattE Integrale \cite{LatInt} for the volume computation.The Condorcet efficiency of plurality voting was first computed by Sch\"urmann in \cite{Sch}, whereas the Condorcet efficiency of the runoff plurality voting was given in \cite{BIS4cand}. According to \cite{GLnew}, it was obtained independently in \cite{BIS4cand} and \cite{Ouafdi}. In \cite[Section 6]{Desc} we additionally discuss the influence of a third ballot on the Condorcet efficiencies of plurality and negative plurality.

Our results for five candidates are listed in Table \ref{computational_results}. The first line contains the probability that the first round winner also wins the runoff. These three computations were done using the full precision mode of Normaliz. The next two lines contain the Condorcet efficiencies, computed the fixed precision mode of Normaliz. For practical reasons we have only included the results rounded to $15$ decimals.
\begin{table}[hbt]
\begin{tabular}{rrrr}
\midrule[1.2pt]
\strut                & \multicolumn{3}{c}{Rule $R$}   \\
\cline{2-4}
\strut                &  \multicolumn{1}{c}{$PR$} & \multicolumn{1}{c}{$NPR$} &  \multicolumn{1}{c}{$BR$} \\
\midrule[1.2pt]
\strut \ttt{RVsRunoff}          & 0.673383666340974 & 0.614598375568014 & 0.769395916647461 \\
\hline
\strut \ttt{CondEffR}           & 0.614270758198443 & 0.509039971570300 & 0.854442922091020 \\
\hline
\strut \ttt{CondEffRRunoff}     & 0.832220522376460 & 0.775488383677566 & 0.991189085613331 \\
\midrule[1.2pt]
\end{tabular}
\vspace*{1ex}
\caption{Probabilities computed by Normaliz} \label{computational_results}
\end{table}

In Table  \ref{Monte_Carlo_results} we reproduce the results for the Condorcet efficiency of all three rules contained in Table 7.6 of \cite{GLnew}, which were obtained using Monte Carlo methods in \cite{LLV}.
\begin{table}[hbt]
\begin{tabular}{rrrr}
\midrule[1.2pt]
\strut                & \multicolumn{3}{c}{Rule $R$}   \\
\cline{2-4}
\strut                &  \multicolumn{1}{c}{$PR$} & \multicolumn{1}{c}{$NPR$} &  \multicolumn{1}{c}{$BR$} \\
\midrule[1.2pt]
\strut \ttt{CondEffR}           & 0.6139 & 0.5090 & 0.8541 \\
\midrule[1.2pt]
\end{tabular}
\vspace*{1ex}
\caption{Results obtained by Monte Carlo, according to \cite{GLnew} and \cite{LLV}} \label{Monte_Carlo_results}
\end{table}
The numbers are relatively close, which confirms the correctness of all algorithms involved. However, at least $14$ decimals printed in Table \ref{computational_results} are exact, while for the numbers printed in Table \ref{Monte_Carlo_results} we have $2$, $4$ and $3$ exact decimals.

\subsection*{Strong Borda paradoxes}

The Borda paradoxes are named after the Chevalier de Borda who studied them in \cite{B}. The \emph{strict Borda paradox} is the event that for a voting profile plurality and majority rank the candidate in opposite order. A less sharp paradox is the \emph{strong Borda paradox}: the plurality winner is the Condorcet loser, and the \emph{reverse strong Borda paradox} occurs if the Condorcet winner finishes last in plurality. These paradoxes can be discussed for all voting schemes for which every profile defines a linear order of the candidates. There is however no point in computing them for negative plurality. As shown in \cite[Section 2.5]{BIS4cand} plurality and negative plurality are dual to each other: the strong Borda paradox and the reverse strong Borda paradox exchange their roles.

For three candidates elections a detailed study of the family of Borda paradoxes \cite{B} is contained in \cite{GL2}, while the case of four candidates is discussed in \cite[Section 2.5]{BIS4cand}. According to \cite{GLnew}, similar results were obtained independently in \cite{Ouafdi}.

For the time being, the computation of the strict Borda paradox in the case of five candidates seems not to be reachable. The strong paradoxes have been computed in the fixed precision mode of Normaliz. The results are rounded to $15$ decimals.

For large numbers of voters the probability of the strong Borda paradox is
$$
B_{\Sg} \approx 0.018125801480904
$$
and the probability of the reverse strong Borda paradox is
$$
B_{\SgRev} \approx 0.019238302806489.
$$

\subsection*{Indifference}\label{Ind_4}
We want to point out that the Normaliz implementation of Lawrence's algorithm does not only yield precise results in five candidates elections, but also extends the range of computations for four candidates considerably by allowing preference rankings with partial indifference that increase the dimension of the related polytopes considerably. We demonstrate this by two examples.

In the examples we allow all possible types of indifference except the equal ranking of all candidates: no indifference, equal ranking of two candidates in three possible positions (top, middle, bottom), two groups of two equally ranked candidates, and equal ranking of three candidates (top and bottom). In total one obtains $74$ rankings. Compared to the $24$ rankings without indifference this is a substantial increase in dimension. We assume that all rankings have the same probability. The authors of \cite{Indiff} allow weights for the types of indifference, for example that the number of voters with a linear order of the candidates is twice the number of voters with indifference. Such weights can easily be realized as a system of homogeneous linear equations in the Normaliz input file.

The first computation is the probability of a Condorcet winner under the Extended Impartial Anonymous Culture (EIAC), as discussed in \cite{Indiff} for $3$ candidates (and varying weights for the different types of indifference). This requires only $3$ inequalities to fix the Condorcet winner, and the computation is very fast. We obtained the value of
$$
0.884041566089553
$$
for the probability of the existence of a Condorcet winner under EIAC (rounded to $15$ decimals).

The second example is the Condorcet efficiency of approval voting. Under this rule one additionally assumes that every voter casts a vote for each candidate on first place in his or her preference ranking. This requires $6$ inequalities, namely $3$ to mark the CW and $3$ to make the same candidate the winner of the approval voting. Consequently the computation time is going up considerably. See the data for \ttt{CondEffAppr 4cand} in Table \ref{calculation_times}. Normaliz obtains
$$
0.695293409282039
$$
as the probability that there exists a CW who finishes first in the approval voting. This yields the Condorcet efficiency of
$$
0.786494024661739
$$
for approval voting (under the assumptions above). The computations were done using the full precision mode of Normaliz.

\subsection*{From three to five candidates}

In Table \ref{3_to_5} we give an overview of the probabilities of voting events for three, four and five candidates as far as we have computed them for five candidates. We use the shorthands $PR$, $NPR$ and $BR$ for the plurality rule, negative plurality rule and Borda rule as introduced above. The remaining abbreviations are self explanatory. For better overview we have rounded all probabilities to 4 decimals.
\begin{table}[hbt]
\begin{tabular}{lrrr}
\midrule[1.2pt]
& 3 cand& 4 cand& 5 cand\\
\midrule[1.2pt]
\strut\ttt{Condorcet Par}             & 0.9375 & 0.8384 & 0.7526\\
\midrule[0.9pt]\strut\ttt{PR vs RunO} & 0.8767 & 0.7545 & 0.6734\\
\hline\strut\ttt{CondEff PR}          & 0.8815 & 0.7426 & 0.6143\\
\hline\strut\ttt{CondEff PR RunO}     & 0.9685 & 0.9117 & 0.8322\\
\midrule[0.9pt]\strut\ttt{NPR vs RunO}& 0.6389 & 0.6227 & 0.6146\\
\hline\strut\ttt{CondEff NPR}         & 0.6296 & 0.5516 & 0.5090\\
\hline\strut\ttt{CondEff NPR RunO}    & 0.9704 & 0.8450 & 0.7755\\
\midrule[0.9pt]\strut\ttt{BR vs RunO} & 0.8750  & 0.8053 & 0.7694\\
\hline\strut\ttt{CondEff BR}          & 0.9111 & 0.8706 & 0.8544\\
\hline\strut\ttt{CondEff BR RunO}     & 1.0000 & 0.9962 & 0.9912\\
\midrule[0.9pt]\strut\ttt{Strong Borda Par} & 0.0296 & 0.0227 & 0.0181\\
\hline\strut\ttt{Strong RevBorda Par} & 0.0315 &0.0238  & 0.0192\\
\midrule[1.2pt]
\end{tabular}
\vspace*{1ex}
\caption{Probabilities of voting events for $3$, $4$ and $5$ candidates}
\label{3_to_5}
\end{table}
One observes that all probabilities are decreasing from three to five candidates. This reflects the increase in the number of configurations defined by the voting profiles. The Condorcet efficiencies and the probabilities of the Borda paradoxes are conditioned on the probabilities of the existence of a Condorcet winner, which itself is decreasing. But this does not compensate the decrease of the absolute probabilities.

In view of our observations above it is justified to formulate

\begin{conjecture}
All series of probabilities associated to voting events in Table \ref{3_to_5} are monotonically decreasing with the number of candidates $n$.
\end{conjecture}

\section*{Condorcet classes}

A voting outcome without ties imposes an asymmetric binary relation on the $n$ candidates that we call a \emph{Condorcet configuration}.
A Condorcet configuration is also called a \emph{dominance relation}, according to \cite{BBH}.
Evidently there are $2^{\binom{n}{2}}$ such configurations.
The permutation group $S_n$ acts on the set of configurations by permuting the candidates. We call the orbits of this action \emph{Condorcet classes}. For $n=4$ the classes and their probabilities are discussed in \cite{BIS4cand}.

From the graph theoretical viewpoint the Condorcet configurations are nothing but simple directed complete graphs with $n$ labeled vertices, i.e., graphs with $n$ labeled vertices without loops, in which each two vertices are connected by a single directed edge. These graphs are also know as \emph{tournament graphs}.
\bigskip

In this section we present the precise probabilities of the Condorcet classes under IAC.
First we make a presentations of the classes, which is needed in order to understand a reduction critical to be made for the computations to be successful.
\bigskip

For $n=5$ these Condorcet configurations fall into $12$ classes under the action of the group $S_5$. There are $6$ classes that have a Condorcet winner (CW) or a Condorcet loser (CL):
\begin{center}
\begin{tabular}{ll}
LinOrd   & CW4cyc\\
CW2nd3cyc\qquad\qquad& 3cyc4thCL\\
CW3cycCL & 4cycCL
\end{tabular}
\end{center}
here ``cyc'' stands for ``cycle''. For example, CW2nd3cyc denotes the class that has a Condorcet winner, a candidate in second position majorizing the remaining three, and the latter are ordered in a $3$-cycle.

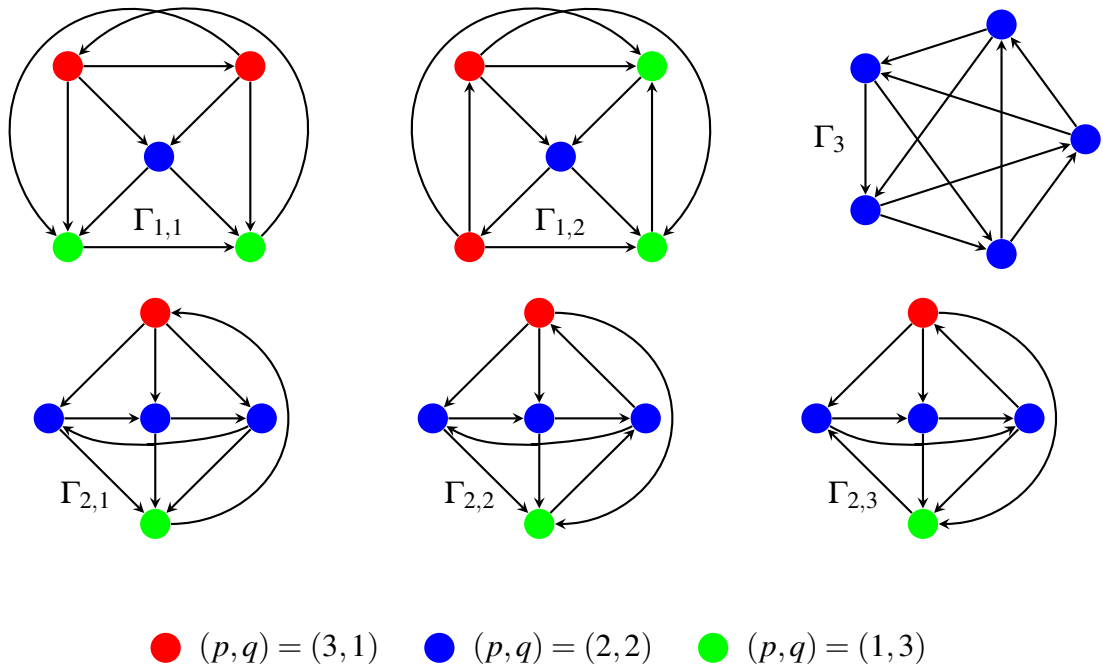
\begin{figure}[hbt]
\begin{center}
\begin{tikzpicture}[squarednode/.style ={rectangle, draw=black!60, very thick, minimum size=5mm}, scale = 0.6]
\tikzstyle{every node} = [circle, fill=gray!30]
\node[fill=white] at (2,0.5){$\Gamma_{1,1}$};
\node[color=red] (a) at (0,4) {};
\node[color=red] (b) at (4,4) {};
\node[color=blue] (c) at (2,2) {};
\node[color=green] (d) at (0,0) {};
\node[color=green] (e) at (4,0) {};
\foreach \from/\to in {a/b, a/c, b/c, c/d, c/e, d/e, a/d, b/e}
\draw [->, >=stealth, thick] (\from) -- (\to);
\draw [->, >=stealth, thick] (e) to [out=45, in=-45] (4.5,4.5) to [out=135, in=45] (a);
\draw [->, >=stealth, thick] (b) to [out=135, in=45] (-0.5,4.5) to [out=225, in=135] (d);
\end{tikzpicture}\qquad
\begin{tikzpicture}[scale = 0.6]
\tikzstyle{every node} = [circle, fill=gray!30]
\node[fill=white] at (2,0.5){$\Gamma_{1,2}$};
\node[color=red] (a) at (0,4) {};
\node[color=green] (b) at (4,4) {};
\node[color=blue] (c) at (2,2) {};
\node[color=red] (d) at (0,0) {};
\node[color=green] (e) at (4,0) {};
\foreach \from/\to in {a/b, a/c, b/c, c/d, c/e, d/e, d/a, e/b}
\draw [->, >=stealth, thick] (\from) -- (\to);
\draw [<-, >=stealth, thick] (e) to [out=45, in=-45] (4.5,4.5) to [out=135, in=45] (a);
\draw [<-, >=stealth, thick] (b) to [out=135, in=45] (-0.5,4.5) to [out=225, in=135] (d);
\end{tikzpicture}\qquad
\begin{tikzpicture}[scale = 1.6]
\tikzstyle{every node} = [circle, fill=blue]
\node[fill=white] at (-1.1,0){$\Gamma_{3}$};
\node (a) at (1,0) {};
\node (b) at (0.309,0.951) {};
\node (c) at (-0.809,0.5878) {};
\node (d) at (-0.809,-0.5878) {};
\node (e) at (0.309,-0.951) {};
\foreach \from/\to in {a/b, b/c, c/d, d/e, e/a, a/c, c/e, e/b, b/d, d/a}
\draw [->, >=stealth, thick] (\from) -- (\to);
\end{tikzpicture}\\[\bigskipamount]
\begin{tikzpicture}[scale = 0.7]
\tikzstyle{every node} = [circle, fill=blue]
\node[fill=white] at (0.7,0.5){$\Gamma_{2,1}$};
\node[color=red] (a) at (2,4) {};
\node (b) at (0,2) {};
\node (c) at (2,2) {};
\node (d) at (4,2) {};
\node[color=green] (e) at (2,0) {};
\foreach \from/\to in {a/b, a/c, a/d, b/c, c/d, b/e, c/e, d/e}
\draw [->, >=stealth, thick] (\from) -- (\to);
\draw [->, >=stealth, thick] (e) to [out=0, in=-90] (4.5,2.0) to [out=90, in=0] (a);
\draw [<-, >=stealth, thick] (b) to [out=-30, in=180] (2.0,1.5) to [out=180, in=-150] (d);
\end{tikzpicture}\qquad\qquad
\begin{tikzpicture}[scale = 0.7]
\tikzstyle{every node} = [circle, fill=blue]
\node[fill=white] at (0.7,0.5){$\Gamma_{2,2}$};
\node[color=red] (a) at (2,4) {};
\node (b) at (0,2) {};
\node (c) at (2,2) {};
\node (d) at (4,2) {};
\node[color=green] (e) at (2,0) {};
\foreach \from/\to in {a/b, a/c, d/a, b/c, c/d, b/e, c/e, e/d}
\draw [->, >=stealth, thick] (\from) -- (\to);
\draw [<-, >=stealth, thick] (e) to [out=0, in=-90] (4.5,2.0) to [out=90, in=0] (a);
\draw [<-, >=stealth, thick] (b) to [out=-30, in=180] (2.0,1.5) to [out=180, in=-150] (d);
\end{tikzpicture}\qquad\qquad
\begin{tikzpicture}[scale = 0.7]
\tikzstyle{every node} = [circle, fill=blue]
\node[fill=white] at (0.7,0.5){$\Gamma_{2,3}$};
\node[color=red] (a) at (2,4) {};
\node (b) at (0,2) {};
\node (c) at (2,2) {};
\node (d) at (4,2) {};
\node[color=green] (e) at (2,0) {};
\foreach \from/\to in {a/b, a/c, d/a, b/c, c/d, e/b, c/e, d/e}
\draw [->, >=stealth, thick] (\from) -- (\to);
\draw [<-, >=stealth, thick] (e) to [out=0, in=-90] (4.5,2.0) to [out=90, in=0] (a);
\draw [->, >=stealth, thick] (b) to [out=-30, in=180] (2.0,1.5) to [out=180, in=-150] (d);
\end{tikzpicture}
\begin{tikzpicture}[squarednode/.style ={rectangle, draw=black!60, very thick, minimum size=5mm}, scale = 0.6]
\tikzstyle{every node} = [circle]
\node[fill=red] (a) at (0,0){};
\node[fill=white] (b)at (2.7,0) {$(p,q)=(3,1)$};
\node[fill=blue] (c) at (6,0){};
\node[fill=white] (d)at (8.7,0){$(p,q)=(2,2)$};
\node[fill=green] (e) at (12,0){};
\node[fill=white] (f) at (14.7,0){$(p,q)=(1,3)$};
\end{tikzpicture}
\end{center}
\caption{The Condorcet classes without a Condorcet winner or loser}\label{CondCl}
\end{figure}

There are $6$ further classes as has been known for a long time. Presumably Davis \cite{Davis} is the oldest source. (For more sources and cardinalities of the set of classes see \cite{oeis}.) The classes can be structured by the \emph{signatures} $(p,q)$ of a candidate in which $p$ counts the candidates majorized by the chosen candidate and $q=n-1-p$ is the number of the candidates majorizing the chosen one. In graph theoretical language, $p$ is the in-degree and $q$ is the out-degree of the chosen node. Without a CW or CL, the signatures $(4,0)$ and $(0,4)$ are excluded. The number of signatures $(2,2)$ must now be odd, and using this observation one easily finds the $6$ classes without a CW or CL. They are named in Figure \ref{CondCl}. In the figure candidates of signature $(3,1)$ are colored red, those of signature $(2,2)$ are blue, and green indicates the signature $(1,3)$.

The cardinalities of all classes and their probabilities (rounded to $6$ decimals) are listed in~
Table~\ref{ClassProb}.

\begin{table}[hbt]
\begin{tabular}{rrr}
\midrule[1.2pt]
\strut   class & \quad$\#$config & \qquad$p$(class)\\
\midrule[1.2pt]
\strut LinOrd& 120&0.533665 \\
\hline
\strut CW2nd3cyc& 40&0.066882\\
\hline
\strut CW3cycCL& 40&0.069984\\
\hline
\strut CW4cyc&120&0.082115 \\
\hline
\strut 3cyc4thCL&40&0.066882\\
\hline
\strut 4cycCL& 120&0.082115\\
\hline
\strut
$\Gamma_{1,1}$ & 120&0.031467\\
\hline
\strut
$\Gamma_{1,2}$& 120&0.032172\\
\hline
\strut
$\Gamma_{2,1}$&40&0.004509\\
\hline
\strut
$\Gamma_{2,2}$&120&0.014644\\
\hline
\strut
$\Gamma_{2,3}$&120&0.014203\\
\hline
\strut
$\Gamma_3$& 24&0.001362\\
\hline
\end{tabular}
\vspace*{1ex}
\caption{Condorcet classes, their cardinalities and probabilities}
\label{ClassProb}
\end{table}

We have computed these probabilities not only for aesthetic reasons: that they sum to~$1$ is an excellent test for the correctness of the algorithm.

For effective computations the following reduction is critical. At first it seems that one must use 10 inequalities representing the relation $>_M$ between the five candidates in addition to the $120$ sign inequalities in order to compute the probability of a single class (or configuration). But computations with $130$ inequalities are currently not reachable on the hardware at our disposal. Some observations help to reduce the number of inequalities, significantly easing the computational load. For example, LinOrd can be (and is) computed with $128$ inequalities if one exploits that it is enough to choose the first two in arbitrary order and the candidate for third place. Once the probability of LinOrd is known, the remaining $5$ classes with a CW or CL can be obtained from the Condorcet paradox ($124$ inequalities), CWand2nd ($126$), CWandCL ($127$) and the symmetry between CW and CL (see \cite{BIS4cand}).

For the other $6$ classes it is best to ``relax'' the direction of some edges and to count which configurations occur if one chooses directions for the relaxed edges. For a proper choice of relaxed edges one gets away with $127$ inequalities for $\Gamma_{1,1}$ and only $126$ or $125$ inequalities for the remaining cases.

It is no surprise that all Condorcet classes have positive probability. In fact, by a theorem of McGarvey \cite{MG} (also see \cite[Theorem 3.1]{BBH}) all Condorcet configurations can be realized by a voting profile. So Proposition \ref{pos_prob} implies positive probability.

The problem of finding the minimal number of voters that are necessary to realize a given Condorcet configuration or even a voting event is largely unknown; see \cite{EM} for an asymptotic lower bound. Some values for four candidates elections have been computed by Normaliz; see \cite[Remark 8]{BIS4cand}.

\section*{Implementations of the Lawrence algorithm and their limitations}\label{algo}

The Lawrence algorithm is based on the fact that a ``signed decomposition'' into simplicies of the polytope in the primal space may be obtained from a ``generic triangulation'' $\Delta$ of its dual cone. For each $\delta\in\Delta$ we get a simplex $R_\delta$ in the primal space and the volume of the polytope in the primal space is the sum of volumes of simplices $R_\delta$ induced by the ``generic triangulation'' with appropriate signs $e(\delta) = \pm 1$. Thus the following formula can be used for computing the volume of $P$:
\begin{equation}
\vol P = \sum_{\delta\in \Delta} e(\delta) \vol R_\delta.\label{iota_vol}
\end{equation}
For mathematical details we refer the reader to Filliman \cite{Fil}. Details of its implementation in Normaliz are described in \cite{BVol}.

In order to compute a ``generic triangulation'', Normaliz, following Lawrence's suggestion, finds a ``generic element'' $\omega$, which in turn induces the ``generic triangulation'' $\Delta = \Delta_\omega$. Since $\omega$ almost inevitably has unpleasantly large coordinates, the induced simplices $R_\delta$ have even worse rational vertices, and their volumes usually are rational numbers with very large numerators and denominators. This extreme arithmetical complexity makes computations with full precision sometimes very difficult on the hardware at our disposal. In the fixed precision mode the volumes $\vol R_\delta$ are computed precisely as rational numbers. But the addition of these numbers may result in gigabytes filling fractions. Therefore in order to make computations feasible the precise rational numbers are truncated to a predetermined set of exact decimal digits, which is typically $100$ digits. Then the error is bounded above by $T\cdot 10^{-100}$ where $T$ is the size of the ``generic triangulation'' (i.e. the total number of simplices).

\begin{remark}\label{fail_vinci}
Before Normaliz, the program vinci \cite{vinci} has provided an implementation of the Lawrence algorithm using  floating point arithmetic. As it is noted by the authors in \cite{practical}, their floating point implementation is numerically unstable. We point out at least one possible reason for this problem, which is indicated by the above discussion.

In any implementation of the Lawrence algorithm the alternating sum \ref{iota_vol} must be evaluated. When using floating point arithmetic for subtracting nearby quantities it is possible that the most significant digits are equal and they will cancel each other. This is a severe limitation of the floating point arithmetic that may lead to a phenomenon known as "catastrophic cancelation".
It is a fact that, because of the relative error involved, the evaluation of a single subtraction in floating point arithmetic could produce completely meaningless digits.

This problem is visible already when computing voting problems with $4$ candidates and only becomes worse for $5$ candidates. Consider the problem of comparing $4$ voting rules for $4$ candidates as it is presented in detail in \cite[Sect.\ 6.1]{Desc}. With its HOT algorithm vinci computes the precise associated Euclidean volume of $1.260510232743\cdot 10^{-25}$. At the same time, a computation with the Lawrence algorithm as it is implemented in vinci provides the erroneous value of $9.287423132835\cdot 10^{-8}$ for the same volume. So is clear that the results provided by the vinci implementation of the Lawrence algorithm may lack any kind of precision, therefore it does not make sense to include in this paper a benchmark of the (different) implementation of the Lawrence algorithm in vinci.
\end{remark}

\begin{remark}\label{polymake} The program polymake \cite{pmk} has also implemented a simplified version of the Lawrence's algorithm.
This implementation is restricted to the "smooth" case. Note that smooth implies "simple", which in turn implies that the dual polytope is "simplicial", so its boundary has a trivial triangulation.
The polytopes that appear in voting theory are not smooth, in fact they are not even simple. Thus the implementation in polymake of the Lawrence algorithm cannot be compared with the Normaliz implementation
for the polytopes presented here.
\end{remark}

\section*{Computational report}\label{computations}

\subsection*{Selected examples}
In order to give the reader an impression of the computational effort, we illustrate it by the data of several selected examples. Except (1) and (2) they are all computations for elections with $5$ candidates:
\begin{enumerate}
\item \ttt{strictBorda 4cand} is the computation of the probability of the strict Borda paradox for elections with $4$ candidates as discussed in \cite{BIS4cand}.

\item \ttt{CondEffAppr 4cand} is the Condorcet efficiency of approval voting for $4$ candidates.

\item \ttt{Condorcet} stands for the existence of a Condorcet winner in elections with $5$ candidates.

\item \ttt{PlurVsRunoff} computes the probability that the plurality winner also wins the runoff.

\item \ttt{CWand2nd} computes the probability that there exists Condorcet winner and a second candidate dominating the remaining three.

\item \ttt{CondEffPlurRunoff} is used to compute the probability that the Condorcet winner exists and finishes at least second in plurality.

\item \ttt{CondEffPlur} computes the probability that the Condorcet winner exists and wins plurality.
\end{enumerate}
In all cases one has to make choices for the candidates that have certain roles in the computation in order to define the polytope for the computation. Table \ref{numerical_data} contains their characteristic combinatorial data.
\begin{table}[hbt]
\begin{tabular}{rrrrr}
\midrule[1.2pt]
\strut                            & dim $C$ & $\#$ inequalities & $\#$ triangulation &  $\#$ generic triang   \\
\midrule[1.2pt]
\strut \ttt{strictBorda 4cand}    &   24    &          33       &            100,738 &         324,862 \\
\hline
\strut \ttt{CondEffAppr 4cand}     &   74    &          80       &              1,620,052 &         30,564,920 \\
\hline
\strut \ttt{Condorcet}            &  120    &         124       &            137,105 &       6,572,904 \\
\hline
\strut \ttt{PlurVsRunoff}         &  120    &         125       &          4,912,369 &      93,749,784 \\
\hline
\strut \ttt{CWand2nd}             &  120    &         126       &         15,529,730 &     608,572,514 \\
\hline
\strut \ttt{CondEffPlurRunoff}    &  120    &         127       &        246,310,369 &   5,456,573,880 \\
\hline
\strut \ttt{CondEffPlur}          &  120    &         128       &      2,388,564,481 & 39,390,184,920 \\
\midrule[1.2pt]
\end{tabular}
\vspace*{1ex}
\caption{Combinatorial data} \label{numerical_data}
\end{table}

\subsection*{Parallelized and distributed volume computations}

The implementation in Normaliz of the Lawrence algorithm consists of $4$ distinct steps that are described in \cite{BVol}.
For effective computations these steps can be separated (and sometimes they must be separated)
and run on different machines.

The computation times in Table \ref{calculation_times} are ``wall clock times'' taken on a Dell R640 system with $1$ TB of RAM and two Intel\texttrademark Xeon\texttrademark Gold 6152 (a total of $44$ cores) using $32$ parallel threads (of the maximum of $88$).
\begin{table}[hbt]
\begin{tabular}{rrrrr}
\midrule[1.2pt]
\strut                 & \multicolumn{1}{c}{RAM} & \multicolumn{3}{c}{time }   \\
\cline{3-5}
\strut                 &  \multicolumn{1}{c}{in GB}&  \multicolumn{1}{c}{stages (1) -- (3)} & \multicolumn{1}{c}{stage (4)} &  \multicolumn{1}{c}{total} \\
\midrule[1.2pt]
\strut \ttt{strictBorda 4cand}     &   0.35 &     1.278 s  &     0.464 s &    1.742 s \\
\hline
\strut \ttt{CondEffAppr 4cand}      &   7.4 &     97.8 s  &     14:31 m &    16:09 m \\
\hline
\strut \ttt{Condorcet}             &   1.67 &    18.0 s  &    52.493 s &     1:10 m \\
\hline
\strut \ttt{PlurVsRunoff}          &  26.2 &     12:40 m  &   1:29:21 h &  1:42:01 s \\
\hline
\strut \ttt{CWand2nd}              &  56.4 &     49:55 m  &  10:21:36 h & 11:11:31 h \\
\hline
\strut \ttt{CondEffPlurRunoff}     & 113 &  13:30:22 h  &   HPC          &  --- \\
\hline
\strut \ttt{CondEffPlur}           & 646    & 125:27:20 h   &   HPC         &  --- \\
\midrule[1.2pt]
\end{tabular}
\vspace*{1ex}
\caption{Memory usage and times for parallelized volume computations} \label{calculation_times}
\end{table}Additional information:
\begin{enumerate}
\item All computations in the table use $64$ bit integers for steps (1)--(3). Even step (4) is done with $64$ bit integers for \ttt{strictBorda 4cand} and \ttt{Condorcet}.

\item The volumes of the first $5$ polytopes were computed with full precision, whereas for \texttt{CondEffPlur} and \texttt{CondEffPlurRunoff} fixed precision was used.

\item The following rule of thumb can be used to estimate the computation time for a smaller number of threads: if one reduces the number of parallel threads from $32$ to $8$, then one should expect the computation time to go up by a factor of $3$. A further reduction to $1$ thread increases it by another factor of $7$.

\item From the selected examples, only \ttt{strictBorda} is computable with the algorithms previously implemented in Normaliz. For this example, the data in Table \ref{calculation_times} may be compared with the data in \cite[Table 2]{Desc} which was recorded on the same system.

\item The data in Table \ref{calculation_times} shows why computations with more than $128$ inequalities are currently not reachable on the hardware at our disposal. Each additional inequality added leads to a significant jump in the required RAM memory and there exists a $1$ TB limit on our system.
\end{enumerate}

Stage (4) of the last two polytopes was computed on a high performance cluster (HPC) because the computation time would become extremely long on the R640, despite of the high degree of internal parallelization. The time for \texttt{CondEffPlurRunoff} would still be acceptable, but \texttt{CondEffPlur} would take several weeks. Instead doing step (4) directly, the result of steps (1)--(3) is written to a series of compressed files on the hard disk. Each of these files contains a certain number of simplices and this number can be chosen by the user, for example $10^6$ simplices. For \texttt{CondEffPlur} we need $12277$ seconds for writing the input files of the distributed computation, and \texttt{CondEffPlurRunoff} needs 528 seconds.

The compressed files are then collected and transferred to the HPC. The Osnabr\"uck HPC has 51 nodes, each equipped with 1 TB of RAM and $2$ AMD Epyc 7742 so that $128$ threads can be run on each node. In our setup each node ran 16 instances of \texttt{chunk} simultaneously and every instance used $8$ threads of OpenMP parallelization. Consequently 816 input files could be processed simultaneously. For a \ttt{CondEffPlur} input file of $10^6$ simplices one needs about $165$ MB of RAM and $3$ hours of computation time. Therefore the volume of \texttt{CondEffPlur} could be computed in $\approx 9$ hours.

Even on a less powerful system it can be advisable to choose this type of approach since one loses only a small amount of data when a system crash should happen and the amount of memory used remains low. Also ``small'' computations can profit from fixed precision. For example, step (4) of \texttt{Condorcet} takes $13.9$ seconds with fixed precision, but $52.5$ seconds with full precision.

\section*{Acknowledgments}

The first author was supported by the DFG (German Research Foundation) grant Br 688/26-1. The second author was partially supported by a grant of Romanian Ministry of Research, Innovation and Digitization,
CNCS/CCCDI - UEFISCDI, project number PN-III-P4-ID-PCE-2020-0878, within PNCDI III.

The high performance cluster of the University of Osnabr\"uck that made the computations possible was financed by the DFG grant 456666331. We cordially thank Lars Knipschild, the administrator of the HPC, for his assistance.

Our thanks also go to Ulrich von der Ohe for his careful reading of the manuscript.

\section*{Author contributions}

All authors have been equally contributing to the paper.

\section*{Data availability}

Input and output files for all computations of this paper can be found at
\begin{center}
\url{https://www.normaliz.uni-osnabrueck.de/documentation/interesting-and-challenging-examples-for-normaliz/}
\end{center}

\section*{Competing interests}

The authors declare no competing interests.

\end{document}